\documentclass[11pt]{article}
\usepackage{amssymb,amsmath,amscd}
\newtheorem{prop}{Proposition}[section]
\newtheorem{lem}{Lemma}[section]
\newtheorem{ex}{Example}[section]
\newtheorem{thm}{Theorem}[section]
\newtheorem{cor}{Corollary}[section]
\newtheorem{defn}{Definition}

\newcommand{\f}[1]{\mathfrak{#1}}

\newcommand{\mb}{\mathbb}
\newcommand{\commentout}[1]{}

\newcommand{\mc}{\mathcal}

\newcommand{\Hom}{\rm Hom\,}
\newcommand{\wint}{\int^{\ast}}

\newcommand{\Di}{{\rm D}_c^{\infty}}
\begin{document}
\title{Generalized Matrix Coefficients for Infinite Dimensional Unitary Representations}
\author{Hongyu He \footnote{AMS Classification 46E, 46F, 22E }\\
Department of Mathematics, Yale University\\
and \\
Department of Mathematics, Louisiana State University \\
email: livingstone@alum.mit.edu\\
}
\date{}
\maketitle
\abstract{Let $(\pi, \mc H)$ be a unitary representation of a Lie group $G$. Classically, matrix coefficients are continuous functions on $G$ attached to a pair of vectors in $\mc H$ and $\mc H^*$. In this note, we generalize the definition of matrix coefficients to a pair of distributions in $(\mc H^{-\infty}, (\mc H^*)^{-\infty})$. Generalized matrix coefficients are in $\mathbf D^{\prime}(G)$, the space of distributions on $G$. By analyzing the structure of generalized matrix coefficients, we prove that, fixing an element in $(\mc H^*)^{-\infty}$, the map $\mc H^{-\infty} \rightarrow \mathbf D^{\prime}(G)$ is continuous. This effectively answers the question about computing generalized matrix coefficients.  For the Heisenberg group, our generalized matrix coefficients can be considered as a generalization of the Fourier-Wigner transform. }

\section{Introduction}
Let $G$ be a Lie group and $H$ a closed subgroup. Let $(\pi, \mc H)$ be a unitary representation of $G$. Throughout this paper, we assume the Hilbert space $\mc H$ to be separable. Let $\mc H^{\infty}$ be the space of smooth vectors equipped with the natural topology defined by semi-norms $\{\| v \|_X= \|\pi(X ) v \| \mid X \in U(\f g) \}$. It is well-known that $\mc H^{\infty}$ is a Frechet space (\cite{wallach} \cite{warner}). Let $(\mc H^*)^{-\infty}$ be the topological dual of $\mc H^{\infty}$. We equip it with
\begin{enumerate}
\item the natural contragredient action $\pi^c$ of $G$;
\item the natural contragredient action $\pi^c$ of the Lie algebra $\f g$;
\item the weak star topology, i.e., the pointwise convergence with respect to each vector in $\mc H^{\infty}$.
\end {enumerate}
 We call vectors in $(\mc H^*)^{-\infty}$ distributions. 
 \begin{enumerate}
 \item Let $\mc H^*$ be the  Hilbert space dual  of $\mc H$, consisting of continuous linear functionals of $\mc H$.  $\mc H^*$ is a subspace of $(\mc H^*)^{-\infty}$ and the contragredient  action $\pi^c$ of $G$ on $(\mc H^*)^{-\infty}$ is the canonical extension of $\pi^c$ on $\mc H^*$. 
 \item Let $(\mc H^*)^{\infty}$ be the Frechet space of smooth vectors in $(\pi^c, \mc H^*)$. Let $\mc H^{-\infty}$ be the topological dual of $(\mc H^*)^{\infty}$. The dual actions $(\pi^c)^c$ of $G$  and $\f g$  on $\mc H^{-\infty}$ are the canonical extensions of $\pi$ on $\mc H^{\infty}$. By abusing notations,  we still denote these actions by $\pi$.
 \end{enumerate}
 An important problem in representation theory is to study $H$ invariant distributions. The structure of the $H$-invariant subspace $((\mc H^*)^{-\infty})^H$ should shed lights on the harmonic analysis on $G/H$.
For example, let $\eta$ be an $H$-invariant distribution in $(\mc H^*)^{-\infty}$. Then there is a natural intertwining operator
$$I_{\eta}: \mc H^{\infty} \rightarrow C^{\infty}(H \backslash G)$$
where $I_{\eta}(v)(g)=\langle \pi(g) v, \eta \rangle = M_{v, \eta}(g)$ and $\langle \, , \, \rangle$ denotes the natural pairing between $\mc H^{\infty}$ and
its dual. In some favorable situations, this intertwining map extends to an isometry from $\mc H$ to $L^2(H \backslash G)$. Then we obtain a subrepresentation of $L^2(H \backslash G)$.  However it remains a difficult problem to  compute $ M_{v, \eta}(g)$ when $v \notin \mc H^{\infty}$. \\
\\
The purpose of this note is to study $\mc M_{\zeta, \eta}$ when both $\eta$ and $\zeta$ are distributions. 
Indeed,  for any
$\eta \in (\mc H^*)^{-\infty}$ and $\zeta \in \mc H^{-\infty}$, we define the notion of generalized matrix coefficient
$\mc M_{\zeta, \eta}$ as a distribution in  $\mathbf D(G)^{\prime}$. Here $\mathbf D(G)^{\prime}$ is the space of distributions equipped with the weak star topology. Recall that for $\zeta \in \mc H$, $\eta \in \mc H^*$, matrix coefficient $ M_{\zeta, \eta}(g)$ is defined to be
$$M_{\zeta, \eta}(g)=\langle \pi(g) \zeta, \eta \rangle=\langle \zeta, \pi^c(g^{-1}) \eta \rangle =\eta(\pi(g) \zeta).$$
In this situation, the distribution $\mc M_{\zeta, \eta}$ coincides with the function $M_{\zeta, \eta}$. Moreover, 
 generalized matrix coefficient satisfies all the covariant properties just like matrix coefficient.
\begin{thm} Let $(\pi, \mc H)$ be a unitary representation of a Lie group $G$. Let $\zeta \in \mc H^{-\infty}$. Then there is a unique continuous $G$-equivariant map $\mc M_{\zeta}: (\mc H^*)^{-\infty} \rightarrow \mathbf D^{\prime}(G)$ such that $\forall \, v \in (\mc H^*)^{\infty}$, $\mc M_{\zeta}(v)$ can be identified with the function
$$\mc M_{\zeta, v}(g)=\langle \pi(g) \zeta, v \rangle =\zeta(\pi^c(g^{-1}) v) .$$
In addition, for any $h \in G$ and $\eta \in (\mc H^*)^{-\infty}$, $\mc M_{\zeta}(\pi^c(h) \eta)= L(h) \mc M_{\zeta}(\eta)$ where $L$ denote the left regular action of $G$ on $\mathbf D(G)^{\prime}$.
If $\pi$ is irreducible and $\zeta \neq 0$, the map $\mc M_{\zeta}$ is injective.
\end{thm}
We shall now make some remarks concerning our theorem.
\begin{enumerate}
\item
The assumption that the representation is unitary is essential in this paper. For $(\pi, V)$ a representation on a topological space, it is impossible to define generalized matrix coefficients on $(V^*, V^*)$ in the most general setting. Nevertheless, if $G$ is semisimple and $\pi$ admissible and finitely generated, there is a well-established theory of Casselman and Wallach which can be used to define generalized matrix coefficients (\cite{ca}\cite{wa}). In this situation, generalized matrix coefficients are more or less independent of the Hilbert structure. 
\item
To answer the question of computing generalized matrix coefficient $\mc M_{\zeta, \eta}$, let $v_i \in (\mc H^*)^{\infty}$ such that $ v_i \rightarrow \eta $ in $(\mc H^*)^{-\infty}$. Then by the result of this paper,
$\mc M_{\zeta, v_i} \rightarrow \mc M_{\zeta, \eta}$ in $\mathbf D(G)^{\prime}$. Observe that $\mc M_{\zeta, v_i}$ is a smooth function on $G$. So Theorem provides an effective way to compute the generalized matrix coefficients. More details of this construction is given in section 8.
\item

One of the most well-known unitary representation is the Schr\"odinger representation of the Heisenberg group $H_n$. In the case, the universal enveloping algebra acts as Weyl algebra. So $\mc H^{\infty}$ is the Schwartz space of rapidly decaying functions $\mc S_n$ and $\mc H^{-\infty}$ is the space of tempered distributions $\mc S^{\prime}_n$.  Generalized matrix coefficient yields a map from
$$\mc M: \mc S_n^{\prime} \times \mc S_n^{\prime} \rightarrow \mathbf D^{\prime}(H_n).$$
Since the center of $H_n$ acts on the representation by scalar, this map can be identified with a map
$$\mc M: \mc S_n^{\prime} \times \mc S_n^{\prime} \rightarrow \mathbf D^{\prime}(\mathbb R^{2n}).$$
This map is classically known as the Fourier-Wigner transform and its image is $\mc S_{2n}^{\prime}$.
\end{enumerate}
\section{Preliminaries}
Let $G$ be a Lie group. Let $(\pi, \mc H)$ be a unitary representation of $G$. Let $\f g$ be the Lie algebra of $G$.
A vector $v \in \mc H$ is said to be smooth if the function 
$$ G \ni g \rightarrow \pi(g)v \in \mc H$$
is a smooth function. Let $\mc H^{\infty}$  be  the space of smooth vectors in $\mc H$. Clearly, $\mc H^{\infty}$ is a linear representation of $G$. Define 
$$ \pi(X) v = \frac{d}{d t}|_{t=0} \pi( \exp t X) v =\lim_{t \rightarrow 0} \frac{\pi(\exp tX) v- v}{t} \qquad ( X \in \f g),$$
where limit is taken under {\it norm convergence}.
Then $\mc H^{\infty}$ becomes a representation of $\f g$, thus a representation of the universal enveloping algebra $U(\f g)$. Fix a basis of $\f g$: $\{ X_1, X_2, \ldots X_l \}$. We adopt the multi-index convention:
$$X^{\alpha}=X_1^{\alpha_1} X_2^{\alpha_2} \ldots X_l^{\alpha_l} \in U(\f g).$$
Equip $\mc H^{\infty}$ with
a countable set of seminorms $\{ \| \pi(X^{\alpha}) v \| \mid \alpha \in \mathbb N^{l}  \}$. Then $\mc H^{\infty}$ is complete under the topology defined by these seminorms. $\mc H^{\infty}$ becomes
a Frechet space. See 4.4.4 of \cite{warner} or 1.6 of \cite{wallach} for details.\\
\\
Let $(\pi^c, \mc H^*)$ be the contragredient unitary representation. We use $\langle \, , \, \rangle$ to denote the pairing between $\mc H$ and $\mc H^*$ or between $\mc H^*$ and $\mc H$. We have
$$ \langle \pi^c(g) u, v \rangle= \langle u, \pi^{-1}(g) v \rangle,  \qquad (u \in \mc H^*, v \in \mc H, g \in G);$$
$$\langle \pi^c(X) u, v \rangle =\langle u, -\pi(X) v \rangle, \qquad (u \in {\mc H^*}^{\infty}, v \in \mc H^{\infty}, X \in \f g).$$
Define a real linear map $i$ from $\mc H^*$ to $\mc H$ by the Riesz representation theorem:
$$\langle u, v \rangle= (v, i(u)), \qquad (u \in \mc H^*, v \in \mc H ).$$
Then $i(\lambda u)= \overline{\lambda} i(u)$. So $i$ defines a conjugate linear isomorphism between $\mc H^*$ and $\mc H$. In addition 
$$(v, i(\pi^c(g)u))=\langle \pi^c(g) u, v \rangle=\langle u, \pi^{-1}(g) v \rangle$$
$$=(\pi^{-1}(g) v, i(u))=(v, \pi(g) i(u)).$$
It follows that $i(\pi^c(g) u)= \pi(g) i(u)$ for any $g \in G$. Consequently, $ i(\pi^c(X) u)=\pi(X) i(u)$ for any $u \in (\mathcal H^*)^{\infty}$ and $X \in \f g$. So $i$ identifies $(\pi, \mc H)$ with $(\pi^c, \mc H^*)$ as real Hilbert space representations. \\
\\
Let $(\mathcal H^*)^{-\infty}$ be the dual space of $\mc H^{\infty}$ equipped with the weak-* topology.  The space $(\mathcal H^*)^{-\infty}$ consists of continuous linear functionals on $\mc H^{\infty}$. The action of $g \in G$ on $\mc H^{\infty}$ induces an action of $g$ on $(\mathcal H^*)^{-\infty}$, still denoted by $\pi^c(g)$,
$$ \langle \pi^c(g) \phi, v \rangle= \langle \phi, \pi^{-1}(g) v \rangle, \qquad (\phi \in (\mathcal H^*)^{-\infty}, v \in \mc H^{\infty}).$$
The action of $X \in \f g$ on $\mc H^{\infty}$ induces an action of $X$ on $(\mathcal H^*)^{-\infty}$, still denoted by  $\pi^c(X)$,
$$\langle \pi^c(X) \phi, v \rangle =\langle \phi, -\pi(X) v \rangle, \qquad (\phi \in (\mathcal H^*)^{-\infty}, v \in \mc H^{\infty}).$$
Given any $X^{\alpha} \in U(\f g)$, define the transpose
$${{}^t X^{\alpha}}=(-1)^{\sum_{i=1}^l \alpha_i} (X_l)^{\alpha_l} X_{l-1}^{\alpha_{l-1}} \ldots X_2^{\alpha_2} X_1^{\alpha_1}.$$
Then 
$$\langle \pi^c(X^{\alpha}) \phi, v \rangle =\langle \phi, \pi({}^t{X^{\alpha}}) v \rangle, \qquad (\phi \in (\mathcal H^*)^{-\infty}, v \in \mc H^{\infty}).$$
We extend the transpose to the universal enveloping algebra $U(\f g)$ by linearality. Transpose defines an anti-automorphism on $U(\f g)$. \\
\\
We retain $D^*$ to denote the conjugate transpose, namely $\overline{{}^t D}$ for any $D \in U(\f g)$. The following proposition is well-known.
\begin{prop}\label{elementary} Let $\phi \in (\mathcal H^*)^{-\infty}$. Then there exists a $u \in \mc H^*$ and $D \in U(\f g)$ such that
$\phi=\pi^c(D) u$ in $(\mathcal H^*)^{-\infty}$.
\end{prop}
Proof: Let $\phi \in (\mathcal H^*)^{-\infty}$. Let $v \in \mc H$. Consider $s_{\phi}(v)=| \langle \phi, v \rangle |$. Then 
$s_{\phi}$ defines a continuous seminorm on $\mc H^{\infty}$. Since the topology of $\mc H^{\infty}$ is generated by the set of seminorms $\{ \| \pi(X^{\alpha}) v \| \mid \alpha \in \mathbb N^l \}$, there exists $D_i \in U(\f g)$, $i \in [1, k]$ such that
$$s_{\phi}(v) \leq \sum_{i=1}^k \| \pi(D_i) v \| = \sum_{i=1}^k ( \pi(D_i) v, \pi(D_i) v )^{\frac{1}{2}}.$$
We have
\begin{equation}
\begin{split}
| \langle \phi, v \rangle | & \leq \sum_{i=1}^k (\pi(D_i^* D_i) v, v)^{\frac{1}{2}} \\
 & \leq \sqrt{k} \, (\sum_{i=1}^k (\pi(D_i^*D_i) v, v))^{\frac{1}{2}} \\
 & \leq \sqrt{k} \, \| \sum_{i=1}^k \pi(D_i^* D_i) v \|^{\frac{1}{2}} \, \| v \|^{\frac{1}{2}} \\
 & \leq \sqrt{k} \, \frac{ \| \sum_{i=1}^k \pi(D_i^* D_i) v \| + \| v \|}{2} \\
 & \leq \sqrt{k} \| \sum_{i=1}^k \pi(D_i^* D_i) v + v \|
\end{split}
\end{equation}
In the last step, we utilize the fact that $(\sum \pi(D_i^* D_i) v, v) \geq 0$. Put
$D=\sqrt{k} (1+\sum_{i=1}^k D_i^* D_i)$. We have $s_{\phi}(v) \leq \| \pi(D) v \|$ for any $v \in \mc H^{\infty}$.\\
\\
Notice that $\|\pi(D) v \| \geq \| v \|$. So $\pi(D): \mc H^{\infty} \rightarrow \mc H^{\infty}$ is injective. 
Now  define a linear functional on $\pi(D) \mc H^{\infty}$ by
$$l(v)=\langle \phi, \pi(D)^{-1} v\rangle \qquad ( v \in \pi(D) \mc H^{\infty}).$$ 
Since $| l(v) |= \| \langle \phi, \pi(D)^{-1} v\rangle \| =s_{\phi}(\pi(D)^{-1} v) \leq \| \pi(D) \pi(D)^{-1} v \|=\| v \|$, $l$ is a bounded linear functional on
$\pi(D) \mc H^{\infty}$. It can be extended to a bounded linear functional on $\mc H$. By Riesz representation theorem, there exists a $u \in \mc H^*$ such that
$$\langle u, \pi(D) v \rangle= l(\pi(D) v)=\langle \phi, \pi(D)^{-1} \pi(D) v \rangle = \langle \phi,  v \rangle, \qquad (v \in \mc H^{\infty}).$$
Notice that the left hand side is $\langle \pi^c({}^t D)u, v \rangle$. We obtain $\phi=\pi^c({}^t D) u$ in $(\mathcal H^*)^{-\infty}$.
$\Box$ \\
\\
So we have proved 
\begin{cor} $(\mathcal H^*)^{-\infty}=\pi^c(U(\f g)) \mc H^*$.
\end{cor}
\begin{ex}
Let $G$ be the Heisenberg group.
Let $\rho$ be the Schr\"odinger representation. The underlying Hilbert space is $L^2(\mathbb R^n)$. The universal enveloping algebra action can be identified with the Weyl algebra. Clearly $\mc H^{\infty}$ is the Schwartz space $\mc S$ and $(\mc H^{*})^{-\infty}$ is $\mc S^{\prime}$, the space of tempered distributions.
Our Corollary simply says that every tempered distribution can be written as $D f$ where $D$ is an algebraic differential operator and $f$ is an $L^2$-function. 
\end{ex}
\section{Matrix Coefficients: Smooth Case}
 Fix $u \in (\mathcal H^*)^{\infty}$. For any $v \in \mc H^{-\infty}$, let 
$$\mc M_{v, u}(g)=\langle \pi(g) v, u \rangle=\langle v, \pi^c(g^{-1}) u \rangle$$
be a matrix coefficient. When $v \in \mc H$, the function $\mc M_{v,u}(g)$ is a smooth function, and the left and right actions of the Lie algebra on $\mc M_{v,u}(g)$ are compatible with the two actions of the Lie algebra on $v$ and $u$ respectively. The purpose of this section is to show that these properties hold when $v \in \mc H^{-\infty}$. \\
\\
Let $L$ be the left regular action of $G$ on $C^{\infty}(G)$ and $R$ be the right regular action of $G$ on $C^{\infty}(G)$. When $v \in \mc H^{-\infty}$ and $h \in G$, we have
$$(R(h)\mc M_{v, u})(g)=\mc M_{v, u}( gh)=\langle \pi(g h) v, u \rangle$$
$$=\langle \pi(g) \pi(h) v, u \rangle=
\mc M_{\pi(h)v, u}(g).$$
$$(L(h)\mc M_{v, u})(g)=\mc M_{v, u}(h^{-1} g)=\langle \pi(h^{-1} g) v, u \rangle$$
$$ =\langle \pi(g) v, \pi^c(h) u \rangle=
\mc M_{v, \pi^c(h) u}(g).$$
So the group actions on the matrix coefficient $\mc M_{v,u}$ is automatically compatible with the groups actions on $v$ and $u$ with $v \in \mc H^{-\infty}$ and $u \in \mc H^*$.\\
\\
The Lie algebra actions require a little bit more caution. 
Recall that for $u \in (\mathcal H^*)^{\infty}$, 
$\pi^c(X) u=\frac{d}{d t}|_{t=0} \pi^c(\exp t X) u$. The limit here is taken with respect to the Hilbert norm.
Since the derivative $\pi^c(X) u$ is defined,
it follows that for any $v \in \mc H$ and $X \in \f g$,
\begin{equation}
\begin{split}
 & (L(X) \mc M_{v, u})(g)=  \frac{d}{d t}|_{t=0} \mc M_{v, u}(\exp (-tX) g) =\frac{d}{ d t}|_{t=0} \langle \pi(g) v, \pi^c(\exp t X) u \rangle \\
& =\langle \pi(g) v, \frac{d}{ d t}|_{t=0} \pi^c(\exp t X) u \rangle= \langle \pi(g) v, \pi^c(X) u \rangle= \mc M_{v, \pi^c(X) u}(g).
\end{split}
\end{equation}
Therefore for any $D \in U(\f g)$, we have $(L(D) \mc M_{v, u})(g)=\mc M_{v, \pi^c(D) u}(g)$ for any $u \in (\mc H^*)^{\infty}$ and $v \in \mc H$. Notice here this simple argument does not work when $v \in \mc H^{-\infty}$. Nevertheless, it is easy to see that $\mc M_{v, u}(g)$ is a smooth function when $u \in (\mc H^*)^{\infty}$ and $v \in \mc H$.\\
\\
Consider the right Lie derivative of $\mc M_{v, u}$. More cautions have to be taken here. For example, it is not obvious that  $(R(D) \mc M_{v, u})(g)=\mc M_{\pi(D) v, u}(g)$ holds for  $v \in \mc H^{-\infty}$ and $u \in (\mc H^*)^{\infty}$.   The main purpose of the following theorem is to extend the covariance of the Lie algebra actions from $v \in \mc H$ to $\phi \in \mc H^{-\infty}$.
\begin{thm} For any $\phi \in \mc H^{-\infty}, u \in (\mathcal H^*)^{\infty}$, the function
$\mc M_{\phi, u}(g)$ is smooth and for any $ D \in U(\f g)$
$$(R(D) \mc M_{\phi, u})(g)=\mc M_{\pi(D) \phi, u}(g), \qquad (L(D) \mc M_{\phi, u})(g)=\mc M_{\phi, \pi^c(D) u}(g).$$
In particular,
$$\frac{d}{ d t}|_{t=0} \langle \pi(\exp t X) \phi, u \rangle= \langle \pi(X) \phi,  u \rangle.$$
\end{thm}
Proof:  For any $X \in \f g$, $u \in (\mc H^*)^{\infty}$, $v \in \mc H$ and $n \in \mathbb Z^+$, we have
\begin{equation}
\begin{split}
\mc M_{\pi(X^n)v, u}(g)= & \langle \pi(g) \pi(X^n) v, u \rangle\\
= & \langle v, \pi^c((-X)^n) \pi^c(g^{-1}) u \rangle \\
= & \langle v, \frac{d^n }{ d t^n}|_{t=0} \pi^c(\exp -t X) \pi^c(g^{-1}) u \rangle \\
= & \frac{ d^n }{ d t^n}|_{t=0} \langle v, \pi^c(\exp (-tX) g^{-1}) u \rangle \\
= & \frac{ d^n }{ d t^n}|_{t=0} \langle \pi(g \exp tX) v, u \rangle \\
= & \frac{ d^n }{ d t^n}|_{t=0} (R(\exp tX) \mc M_{v, u})(g) \\
= & R(X^n) \mc M_{v, u}(g) \\
\end{split}
\end{equation}
\noindent
\noindent
\noindent
\noindent
By linearality, $\langle \pi(g) \pi(D) v, u \rangle=R(D) \langle \pi(g) v, u \rangle$ for any $D \in U(\f g)$. Notice that the latter is always smooth. So $\mc M_{\pi(D) v, u}(g)$ is a smooth function. By Proposition \ref{elementary}, any $\phi \in \mc H^{-\infty}$ can be written as $\pi(D_1) v$ for some $v \in \mc H$. We see that $\mc M_{\phi, u}(g)$ is a smooth function. 
In addition,  we have
$$R(D) \mc M_{\phi, u}(g)= R(D) R(D_1) \mc M_{v, u}(g)=\mc M_{\pi(D D_1) v, u}(g)=\mc M_{\pi(D) \phi, u}(g).$$
Similarly, we have
\begin{equation}
\begin{split}
& L(D) \mc M_{\phi, u}(g) =L(D) R(D_1) \mc M_{v, u}(g)=R(D_1) L(D) \mc M_{v,u}(g) \\
= & R(D_1) \mc M_{v, \pi^c(D) u}(g) 
 =  \mc M_{\pi(D_1) v, \pi^c(D) u}(g)=\mc M_{\phi, \pi^c(D) u}(g).
 \end{split}
 \end{equation}
\noindent
It follows that
$$\frac{d}{ d t}|_{t=0} \langle \pi(\exp t X) \phi, u \rangle=(R(X) \mc M_{\phi, u})(e)= \mc M_{\pi(X) \phi, u}(e)= \langle \pi(X) \phi,  u \rangle.$$ $\Box$

\section{Weak Integral of distributions}
Let $(X, \mu)$ be a measure space, $V$ a topological vector space and $\Phi: X \rightarrow V$. We would like to define $\int_X \Phi(x) d \mu(x)$ in a proper sense. There are various ways this can be accomplished in different settings. We adopt the following natural definition. See \cite{folland}, \cite{rudin}.
\begin{defn} Let $(X, \mu)$ be a measure space and $V$ a locally convex topological vector space. Let $V^*$ be the dual of $V$ equipped with the weak-* topology. Let $\Phi: X \rightarrow V^*$ be a map such that
$\int_X \langle \Phi(x), v\rangle d \mu(x)$ converges absolutely for each $ v \in V$. Then we define
$\wint_X  \Phi(x) d \mu(x)$ as a vector in  $\Hom(V, \mathbb C)$:
$$[\wint_X \Phi(x) d \mu(x) ](v)=\int_X \langle \Phi(x), v\rangle d \mu(x).$$
If {\it in addition} $\wint_X \Phi(x) d \mu(x) \in \Hom(V, \mathbb C)$ is in $V^*$, we say that $\wint_X \Phi(x) d \mu(x)$ converges in $V^*$ or converges weakly.
\end{defn}
\begin{ex} Let $F$ be a tempered distribution on $\mathbb R$. Suppose that its Fourier transform $\hat F$ is a locally integrable function. Since $\hat F \in \mc S^{\prime}(\mathbb R)$, $\hat F$ will be of at most polynomial growth. Then we will always have
$$F= \int^* \hat F(\xi) \exp 2 \pi i x \xi d \xi \qquad {\rm in} \hspace{.05 in} \mc S^{\prime}(\mb R),$$
even though $\int \hat F(\xi) \exp 2 \pi i x \xi d \xi$ may not converge for any $x$. It follows that
$$F= \int^* \hat F(\xi) \exp 2 \pi i x \xi d \xi \qquad {\rm in} \hspace{.05 in} \mathbf D^{\prime}(\mb R).$$
\end{ex}
\begin{ex} 
Let $\mathbb T=\mathbb R/\mathbb Z$ be the 1 dimensional torus. Let $F$ be a distribution in $\mathbf D^{\prime}(\mathbb T)$. Then $\hat F$ as a function on $\mathbb Z$ is of at most polynomial growth. Even though $\sum \hat F(n) \exp 2 \pi i n t$ may not be summable, we will always have
$$F= \sum^* \hat F(n) \exp 2 \pi i n t \qquad {\rm in} \hspace{.05 in} \mathbf D^{\prime}(\mathbb T).$$
\end{ex}
If $V$ is a Hilbert space, $\int_X \Phi(x) d x$ may be defined directly in the Hilbert space, under certain continuous or measurable condition. More generally, when $V$ is a Banach space, there is a well-defined notion of Bochner integral when
$\Phi$ is Bochner measurable and $\int_X \| \Phi(x) \| d x < \infty$. Throughout this paper, $\Phi(x)$ will always be continuous.  So $\wint_X \Phi(x) d x$ defined in the weak sense coincides with
the Bochner integral $\int_X \Phi(x) d x$ when $V$ is a Hilbert space. 
\section{Main Proposition}
Let $G$ be a  Lie group with the {\it left Haar measure} $dg$. Let $(\pi, V)$ be a continuous representation of $G$ on a {\it complete} locally convex topological vector space $V$. Let $f \in C_c^{\infty}(G)$. For any $v \in V$, define
$\pi(f)(v)=\int_G f(g) \pi(g) v d g$. Then $\pi(f)v$ is well-defined. Notice that $V$ is complete implies that $\pi(f)(v)$ converges in seminorm. In addition, the operator $\pi(f): V \rightarrow V$ is continuous. See \cite{warner}. \\
\\
Let $V^*$ be the topological dual space of $V$ equipped with the weak star topology. It may not be complete. Let $(\pi^c, V^*)$ be the dual representation, not necessarily continuous. Nevertheless the action of 
$C_c^{\infty}(G)$ on $V$ still induces an action of $C_c^{\infty}(G)$ on $V^*$, namely
$$\pi^c(f) \phi= \wint_G f(g) \pi^c(g) \phi d g.$$
To see this, observe
$$\langle \pi^c(f) \phi, v \rangle=\int f(g) \langle \phi, \pi(g^{-1}) v \rangle d g =\langle \phi, \int f(g) \pi(g^{-1}) v d g \rangle,$$
$$(\phi \in V^*, v \in V).$$
So $\pi^c(f)  \phi$ is defined in $\Hom(V, \mathbb C)$. Since the operator $v \rightarrow \int f(g) \pi(g^{-1}) v dg$ is continuous  on $V$, 
$\pi^c(f) \phi \in V^*$ and $\pi^c(f): V^* \rightarrow V^*$ is continuous.  See Prop 19.5 and its corollary in \cite{tr}.    \\
\\
In practice, we can use 
$$\langle \pi^c(f) \phi, v \rangle =\langle \phi, \int f(g) \pi(g^{-1}) v d g \rangle, \qquad (\phi \in V^*, v \in V)$$
as the definition of $\pi^c(f) \phi$. 
The main result of this section is
\begin{prop}\label{main1}
Let $(\pi, \mc H)$ be a unitary representation of a Lie group $G$. Let $C_c^{\infty}(G)$ act on $\mc H^{-\infty}$. Let $\phi \in \mc H^{-\infty}$ and $f \in C_c^{\infty}(G)$. 
\begin{enumerate}
\item For any $D \in U(\f g)$, $\pi(D) \pi(f) \phi=\pi(L(D) f) \phi$;
\item for any $D \in U(\f g)$, there is a well-defined anti-homomorphism $\mc A: U(\f g) \rightarrow U(\f g)$ such that $\pi(f) \pi(D) \phi= \pi(R(\mc A( D))f) \phi$ (when $G$ is unimodular, $\mc A(D)={}^t D$);
\item  $\pi(f) \phi \in \mc H^{\infty}$;
\item the map $\mc C_{\phi}: C_c^{\infty}(G) \rightarrow \mc H^{\infty}$ defined by
$$\mc C_{\phi}(f)=\pi(f)(\phi)$$
is continuous.
\end{enumerate}
\end{prop}
This proposition says that convolution with $C_c^{\infty}(G)$ smoothens distributions. The last statement says that $\mc C_{\phi}$ can be regarded as $\mc H^{\infty}$-valued distribution. Before we start our proof, let $\Delta(h)$ be the modular function, namely $ d (g h)= \Delta(h) d g$ for any fixed $h \in G$. Let $\delta$ be the derivative of $\Delta(h)$ at the identity, namely $\delta(X)=\frac{d}{d t}|_{t=0} \Delta(\exp t X)$. Both $\Delta$ and $\delta$ can be computed explicitly for Lie groups. For each $X \in \f g$, we define $\mc A(X)=-X^t-\delta(X)$. Extend $\mc A$ to an anti-automorphism of $U(\f g)$. \\
\\
We begin our proof with the following lemma. These lemma are all well-known. 
\begin{lem}\label{Gar} Let $G$ be a Lie group. Let $\phi \in \mc H^{-\infty}$, $u \in \mc H$ and  $f \in C_c^{\infty}(G)$ . Then 
$$\pi(h) \pi(f) \phi= \pi(L(h) f) \phi, \qquad (h \in G);$$
$$\pi(D) \pi(f) u=\pi(L(D) f) u \qquad (D \in U(\f g)).$$
\end{lem}
Proof: We have
\begin{equation}
\begin{split}
\pi(h)\pi(f) \phi= & \pi(h) (\int^* f(g) \pi(g) \phi d g) \\
= & \int^* f(g) \pi(h) \pi(g) \phi d g \\
= & \int^* f(g) \pi(h g) \phi d g \\
= & \int^* f(h^{-1} g) \pi(g) \phi d g \\
= & \int^* (L(h)f)(g) \pi(g) \phi dg \\
=& \pi(L(h)f) \phi
\end{split}
\end{equation}
All these equations hold when they are evaluated at every smooth vector.
The first statement is proved. The second statement implies the existence of Garding space. Its proof can be found in many textbooks. See \cite{knapp}. $\Box$
\commentout{
Proof: We have for any $v \in {\mc H^{*}}^{\infty}$,
$$\langle \pi(h) \pi(f) \phi, v \rangle=\int_G f(g) \langle \pi(h g) \phi, v \rangle d g= \int_G f(h^{-1} g) \langle \pi(g) \phi, v \rangle= \langle \pi(L(h) f) \phi, v \rangle.$$
By definition, $\pi(h) \pi(f) \phi= \pi(L(h) f) \phi, \qquad (h \in G)$. \\
\\
Suppose that $X \in \f g$, $u \in \mc H$ and $f \in C_c^{\infty}(G)$. We have 
\begin{equation}
\begin{split}
\pi(X) \pi( f) u = & \pi(X) \int_G f(g) \pi(g) u d g \\
= & \frac{d }{d t} |_{t=0} \pi(\exp t X) \int_G f(g) \pi(g) u d g \\
= & \frac{d}{d t}|_{t=0}  \int_G f(g)  \pi(\exp tX g) u d g \\
= & \frac{d} { d t} |_{t=0} \int_G f(\exp (-tX) g) \pi(g) u d g \\
= &  \int_G  \frac{d} { d t} |_{t=0} f(\exp (-tX) g) \pi(g) u d g \\
= & \int_G (L(X) f)(g) \pi(g) u d g \\
= & \pi(L(X) f ) u
\end{split}
\end{equation}
where $L$ denotes the left regular action. It follows that $\pi(X^{\alpha}) \pi(f) u= \pi(L(X^{\alpha}) f) u$ for
$u \in \mc H$ and $f \in C_c^{\infty}(G)$. Hence $\pi(D) \pi(f) u=\pi(L(D) f) u$ for any $D \in U(\f g)$.
$\Box$}
\begin{lem}\label{rightaction} Let $G$ be a  Lie group. Let $u \in \mc H^{\infty}$, $f \in C_c^{\infty}(G)$ and $\phi \in \mc H^{-\infty}$. Then
$$\pi(f) \pi(h) \phi= \Delta(h^{-1}) \pi( R(h^{-1}) f) \phi, \qquad ( h \in G),$$
$$\pi(f) \pi(D) u=\pi(R(\mc A( D) f)) u, \qquad (D \in U(\f g)).$$
\end{lem}
Proof: We have for any $v \in (\mathcal H^*)^{\infty}$,
$$\langle \pi(f) \pi(h) \phi, v \rangle= \langle \int^* f(g) \pi(g) \pi(h) \phi d g, v \rangle =\int f(g) \langle \pi(g h) \phi, v \rangle d g$$
$$= \int f(g h^{-1}) \langle \pi(g) \phi, v \rangle \Delta(h^{-1}) d g= \Delta(h^{-1}) \langle \pi(R(h^{-1}) f) \phi, v \rangle.$$
Hence $\pi(f) \pi(h) \phi= \Delta(h^{-1}) \pi( R(h^{-1}) f) \phi$. \\
\\
Suppose that $X \in \f g$, $u \in \mc H^{\infty}$ and $f \in C_c^{\infty}(G)$. Then we have
\begin{equation}
\begin{split}
 & \pi(f)\pi(X) u  \\
= & \int f(g) \pi(g) \frac{d} { d t}|_{t=0} \pi(\exp t X) u d g \\
 = & \frac{d }{ d t}|_{t=0} \int f(g) \pi( g \exp tX) u d g \\
 = & \frac{d }{ d t}|_{t=0} \Delta(\exp(-t X)) \int f(g \exp (-t X)) \pi(g) u d g \\
 = & \frac{d} { dt}|_{t=0} (\int  f(g \exp (- t X)) \pi(g) u dg )-\delta(X) \int f(g) \pi(g) u d g  \\
 = & \int \frac{d} { dt}|_{t=0} f(g \exp (- t X)) \pi(g) u dg -\delta(X) \int f(g) \pi(g) u d g  \\
 = & \int R(-X) f(g) \pi(g) u d g -\delta(X) \pi(f) u\\
 = & \pi(R(\mc A(X)) f) u
\end{split}
\end{equation}
It follows that for any $D \in U(\f g)$, $\pi(f) \pi(D) u=\pi(R(\mc A( D) )f) u.$
$\Box$ \\
\\
Proof of Proposition \ref{main1}: Let $\phi \in \mc H^{-\infty}$ and $f \in C_c^{\infty}(G)$. Let $v \in (\mathcal H^*)^{\infty}$. Let $X \in \f g$ and $n \in \mathbb N$. Then we have
\begin{equation}
\begin{split}
\langle \pi(X^n) \pi(f) \phi, v \rangle = & \langle \pi(f) \phi, \pi^c((-X)^n) v \rangle \\
=& \langle \int f(g) \pi(g) \phi d g, \pi^c((-X)^n) v \rangle \\
=& \langle \phi, \int f(g) \pi^c(g^{-1}) \pi^c((-X)^n) v d g \rangle \\
=&  \langle \phi, \int f(g) \pi^c(g^{-1}) \frac{d^n}{d t^n}|_{t=0} \pi^c(\exp(-t X)) v d g \rangle  \\
=& \langle \phi, \frac{d^n}{d t^n}|_{t=0} \int f(g) \pi^c(g^{-1} \exp(-t X)) v d g \rangle  \\
=&   \langle \phi, \frac{d^n}{d t^n}|_{t=0} \int f(\exp (-tX) g) \pi^c(g^{-1}) v d g \rangle  \\
=&  \langle \phi, \int \frac{d^n}{d t^n}|_{t=0} f(\exp (-tX)) g) \pi^c(g^{-1}) v d g \rangle \\
=& \langle \phi, \int (L(X^n) f)(g) \pi^c(g^{-1}) v d g \rangle \\
=& \langle \pi(L(X^n) f) \phi, v \rangle 
\end{split}
\end{equation}
Since $U(\f g)$ is generated by $\f g$, by the process of symmetrization, it is spanned by $\{ X^n \mid X \in \f g, n \in \mathbb N \}$. So Prop. \ref{main1} $(1)$ is proved. \\
\\
To prove Prop. \ref{main1} $(2)$, observe that
\begin{equation}
\begin{split}
 &\langle \pi(f)  \pi(X) \phi, v \rangle \\
= & \langle \pi(X) \phi, \int f(g) \pi^c(g^{-1}) v d g \rangle \\
= & \langle \phi, \pi^c(-X) \int f(g) \pi^c(g^{-1}) v d g \rangle \\
= &  \langle \phi, \frac{d }{ d t}|_{t=0} \pi^c(\exp (-tX)) [ \int f(g) \pi^c(g^{-1}) v d g ] \rangle  \\
=& \langle  \phi, \frac{d }{ d t}|_{t=0}  [ \int f(g) \pi^c(\exp (-tX) g^{-1}) v d g ] \rangle \\
= & \langle  \phi, \frac{d }{ d t}|_{t=0}  [ \int f(g \exp (-t X)) \pi^c( g^{-1}) v \Delta(\exp(-tX)) d g ]  \rangle \\
=& \langle \phi, -\delta(X) \int f(g) \pi^c( g^{-1}) v d g \rangle + \langle \phi,  \frac{d }{ d t}|_{t=0} \int  f(g \exp (-t X)) \pi^c( g^{-1}) v d g \rangle \\
=& \langle \phi, -\delta(X) \int f(g) \pi^c( g^{-1}) v d g \rangle + \langle \phi, \int \frac{d }{ d t}|_{t=0} f(g \exp (-t X)) \pi^c( g^{-1}) v d g \rangle \\
=& \langle \pi(-\delta(X) f) \phi, v \rangle - \langle \phi, \int R(X)f(g) \pi^c( g^{-1}) v d g \rangle \\
= & \langle \pi(-\delta(X)-R(X)f) \phi, v \rangle
\end{split}
\end{equation}
By induction, $\pi(f)  \pi(X^n) \phi=\pi( R(\mc A(X)^n) f) \phi$. Prop. \ref{main1} $(2)$ follows immediately.\\
\\ 
By Prop. \ref{elementary}, let $\phi=\pi(D) u$ for some $D \in U(\f g)$ and $u \in \mc H$. Then $\pi(f) \phi=\pi(f) \pi(D) u=\pi(R(\mc A( D)) f) u$. Notice that $R(\mc A(D)) f \in C_c^{\infty}(G)$. By the Theorem of Garding, $\pi(R(\mc A( D)) f) u \in \mc H^{\infty}$. See Lemma \ref{Gar}. So $\pi(f) \phi \in \mc H^{\infty}$. \\
\\
Let $\mc C_{\phi}(f)=\pi(f) \phi$. 
  To show that $\mc C_{\phi}: C_c^{\infty}(G) \rightarrow \mc H^{\infty}$ is continuous, it suffices to show that for any compact $K$ in $G$ and any sequence $f_i \rightarrow 0$ in $C_c^{\infty}(K)$, $\mc C_{\phi}(f_i) \rightarrow 0$ in the Frechet space $\mc H^{\infty}$.
Notice that $f_i \rightarrow 0$ in $C_c^{\infty}(K)$ means that all derivatives $\| L(X^{\alpha}) R(X^{\beta}) f_i \|_{sup} \rightarrow 0$. Here $\| * \|_{sup}$ denote the sup-norm. Let $\phi=\pi(D_1) u$ for some $D_1 \in U(\f g)$ and $u \in \mc H$. For any $D \in U(\f g)$, we have
$$
\| \pi(D) \pi(f_i) \phi \|= \| \pi( L(D) f_i) \pi(D_1) u \|=\|\pi(L(D) R(\mc A( D_1)) f_i) u \| $$
$$ \leq \|L(D) R(\mc A( D_1)) f_i\|_{sup} |K| \|u \| \rightarrow 0.
$$
where $|K|$ denotes the measure of $K$. Therefore $\mc C_{\phi}: C_c^{\infty}(G) \rightarrow \mc H^{\infty}$ is continuous. $\Box$
\section{Generalized Matrix Coefficients}
Let $\phi \in \mc H^{-\infty}$ and $\psi \in (\mc {H^*})^{-\infty}$. Since the operator $\mc C_{\phi}: C_c^{\infty}(G) \rightarrow \mc H^{\infty}$ is continuous, we immediately obtain a continuous dual operator for the dual spaces:
$$\mc M_{\phi}: (\mathcal H^*)^{-\infty} \rightarrow \mathbf D^{\prime}(G).$$
More precisely, we have 
$$ \langle  \mc M_{\phi}(\psi), f \rangle=\langle  \psi, \mc C_{\phi}(f)\rangle=\langle  \psi, \pi(f) \phi \rangle.$$
See Prop 19.5 and its corollary in \cite{tr} for the theory of dual operators on dual topological vector spaces.
\begin{defn} Let $\phi \in \mc H^{-\infty}$ and $\psi \in (\mc {H^*})^{-\infty}$. We define
$\mc M_{\phi, \psi}$ to be a distribution on $G$:
$$\langle \mc M_{\phi, \psi}, f \rangle=\langle \pi(f) \phi, \psi \rangle=\langle \mc M_{\phi}(\psi), f \rangle, \qquad (f \in C_c^{\infty}(G)).$$
\end{defn}
Now we have a map
$$\mc M: \mc H^{-\infty} \otimes (\mathcal H^*)^{-\infty} \rightarrow \mathbf D^{\prime}(G).$$
\begin{lem} If $\phi \in \mc H$ and $\psi \in \mc H^*$, then
$\mc M_{\phi, \psi}$ can be identified with the function 
$$M_{\phi, \psi}(g)= \langle \pi(g) \phi, \psi \rangle.$$
\end{lem}
Proof: For any $f \in C_c^{\infty}(G)$, we have  $$\langle \mc M_{\phi, \psi}, f \rangle=\langle \pi(f) \phi, \psi \rangle= \int f(g) \langle \pi(g) \phi, \psi \rangle d g .$$
The lemma follows immediately. $\Box$\\
\\
So $\mc M_{\phi, \psi}(g)$ generalizes the notion of matrix coefficients. Similar statement holds for $\phi \in \mc H^{-\infty}$ and $\psi \in (\mc H^*)^{\infty}$.  We summarize our discussion in the following theorem.
\begin{thm}\label{5.1} For any $\phi \in \mc H^{-\infty}$, $\mc M_{\phi}$ is a continuous map from
$(\mc H^*)^{-\infty}$ to $\mathbf D^{\prime}(G)$ that coincides with the classical definition of matrix coefficients when restricted to $(\mc H^*)^{\infty}$. 
\end{thm}
Let $L$ be the left action of the  Lie group $G$ and $U(\f g)$ on $\mathbf D^{\prime}(G)$ and $R$ be the right action on $\mathbf D^{\prime}(G)$. The following lemma can be established easily the same way as Lemma \ref{Gar} and \ref{rightaction}.
\begin{lem}\label{rightaction1} Let $F \in \mathbf D^{\prime}(G), f \in C_c^{\infty}(G), h \in G, D \in U(\f g)$. We have
$$\langle L(h) F, L(h) f \rangle=\langle F, f \rangle= \Delta(h) \langle R(h) F, R(h) f \rangle;$$
$$\langle L(D) F, f \rangle=\langle F, L({}^t D )f \rangle, \qquad \langle R(D) F, f \rangle= \langle F, R(\mc A( D)) f \rangle.$$
\end{lem}
By combining this lemma with Prop. \ref{main1}, we obtain
\begin{prop}\label{main2} Let $h \in G$, $D \in U(\f g)$, $\phi \in \mc H^{-\infty}$ and $\psi \in (\mathcal H^*)^{-\infty}$. Then generalized matrix coefficients have the following properties
\begin{enumerate}
\item $\mc M_{\pi(h) \phi, \psi}=R(h) \mc M_{\phi, \psi};$
\item $\mc M_{\phi, \pi^c(h) \psi}= L(h) \mc M_{\phi, \psi};$
\item $\mc M_{\phi, \pi^c(D) \psi}= L(D) \mc M_{\phi, \psi};$
\item $\mc M_{\pi(D) \phi, \psi}=R(D) \mc M_{\phi, \psi}.$
\end{enumerate}
\end{prop}
Proof: By Lemma \ref{Gar}, \ref{rightaction} and \ref{rightaction1}, we have
$$\langle \mc M_{\pi(h) \phi, \psi}, f \rangle=\langle \pi(f) \pi(h) \phi, \psi \rangle=\langle \Delta(h^{-1}) \pi(R(h^{-1}) f) \phi, \psi \rangle $$
$$=\Delta(h^{-1}) \langle \mc M_{\phi, \psi}, R(h^{-1}) f \rangle=\langle R(h) \mc M_{\phi, \psi}, f \rangle.$$
The other statements can be established similarly. $\Box$ \\
\\
Let ${\mc M}_{\pi}$ be the linear span of generalized matrix coefficients of a unitary representation $(\pi, \mc H)$.
Let $M_{\pi}$ be the linear span of matrix coefficients of $(\pi, \mc H)$. The following theorem is the direct consequence of Prop. \ref{main1} and Prop. \ref{main2}.
\begin{thm} We have
${\mc M}_{\pi}=L(U(\f g)) R(U(\f g)) M_{\pi}.$
\end{thm}
\section{Basic Properties of Generalized Matrix Coefficients}
\subsection{Semi-invariant Distributions}
Let $H$ be a subgroup of $G$. Then $H$ acts on $\mc H^{-\infty}$. Let $\chi$ be a character of $H$, namely a one-dimensional representation of $H$. Let
$$(\mc H^{-\infty})^{H, \chi}= \{ \phi \in \mc H^{-\infty} \mid \pi(h) \phi=\chi(h) \phi,  \,\, \forall \, h \in H \}.$$
 $(\mc H^{-\infty})^{H, \chi}$ is called the space of {\it semi-invariant distributions}, with respect to $(H, \chi)$. As a consequence of Prop. \ref{main2}, we have

 \begin{thm}
 Let $\phi \in (\mc H^{-\infty})^{H, \chi}$ and $\psi \in (\mathcal H^*)^{-\infty}$. Then
 $R(h) \mc M_{\phi, \psi}= \chi(h) \mc M_{\phi, \psi}.$
 \end{thm}
 In particular, existence of a nonzero semi-invariant distribution $\eta \in (\mc H^{-\infty})^{H, \chi}$ implies  there is an intertwining operator
 $$\mc I^{\infty}_{\eta} : (\mathcal H^*)^{\infty} \rightarrow  C^{\infty}(G/H, \chi).$$
{\it By our results, this intertwining map can be extended to an intertwining operator from $(\mathcal H^*)^{-\infty} $ to $\mathbf D^{\prime}(G)^{R(H), \chi}$}.
 \subsection{Orthogonality}
 Let $\phi, \psi$ be two vectors in $\mc H^{-\infty}$ and $(\mathcal H^*)^{-\infty}$. Obviously $\psi$ can be identified with a vector in $\mc H^{-\infty}$. We say that $\phi \perp_G \psi$ if $\mc M_{\phi, \psi}(g)=0$ in $\mathbf D^{\prime}(G)$. By definition, we have
 \begin{lem}  $\phi \perp_G \psi$ if and only if $\langle \pi(f) \phi, \psi \rangle=0$ for all $f \in C_c^{\infty}(G)$.
 \end{lem}
 We call the closure of $\{ \pi(f) \phi \mid f \in C_c^{\infty}(G) \}$ {\it the Hilbert space generated by $\phi$ }(with respect to $G$). Then
 $\phi \perp_G \psi$ if and only if the Hilbert spaces generated by $\phi$ and $\psi$ are perpendicular. In some cases,  studying a vector in $\mc H^{-\infty}$ may shed light on the structure of the Hilbert space it generates.
 \subsection{Subrepresentations}
 Let $\mc V$ be a subrepresentation of $\mc H$. Let $\mc W$ be its orthogonal complement. Then $\mc V^{\infty} \subseteq \mc H^{\infty}$. In fact we have
 \begin{thm}
 $\mc H^{\infty}=\mc V^{\infty} \oplus  \mc W^{\infty}$ and $\mc H^{-\infty} \cong \mc V^{-\infty} \oplus \mc W^{-\infty}$.
 \end{thm}
 Proof: Obviously $\mc V^{\infty} \oplus \mc W^{\infty} \subseteq \mc H^{\infty}$. To show that $\mc V^{\infty} \oplus \mc W^{\infty} \supseteq \mc H^{\infty}$, let $P_{\mc V}$ be the projection operator onto $\mc V$. Let $v \in \mc H^{\infty}$. Then 
$$\lim_{t \rightarrow 0} \frac{\|\pi(\exp t X) v - v- t \pi(X) v\|}{ t}=0.$$
It then follows that
$$\lim_{t \rightarrow 0} \frac{\|\pi(\exp t X) P_{\mc V} v - P_{\mc V} v- t P_{\mc V} \pi(X) v\|}{ t}=0.$$
So $\pi(X) P_{\mc V} v= P_{\mc V} \pi(X) v$. Similarly, $\pi(X) P_{\mc W} v= P_{\mc W} \pi(X) v$. By induction, $P_{\mc V} v \in \mc V^{\infty}$ and $P_{\mc W} v \in \mc W^{\infty}$. Since $v=P_{\mc V} v + P_{\mc W} v$, 
we have $\mc V^{\infty} \oplus \mc W^{\infty} \supseteq \mc H^{\infty}$. So  $\mc H^{\infty}=\mc V^{\infty} \oplus  \mc W^{\infty}$. \\
\\
Essentially, we have $P_{\mc V}$ and $P_{\mc W}$ commutes with the action of $U(\f g)$ on $\mc H^{\infty}$. The same is then true on $\mc H^{-\infty}$. Recall that $\mc H^{-\infty}=\pi(U(\f g)) \mc H$. It follows easily that any $\phi \in \mc H^{-\infty}$ can be written as a sum of
$\pi(D) P_{\mc V} v + \pi(D) P_{\mc W} v$. So we obtain a map  $\mc H^{-\infty} \rightarrow  \mc V^{-\infty} \oplus \mc W^{-\infty}$.
It is easy to see that this map is surjective and injective. $\Box$
 \subsection{Injectivity}
 Generally speaking, given $\eta \in \mc H^{-\infty}$, the map $\mc M_{\eta}: (\mc H^*)^{-\infty} \rightarrow \mathbf D^{\prime}(G)$ may not be injective.
 However, if $(\pi, \mc H)$ is irreducible, we have the following
 \begin{thm} Suppose $(\pi, \mc H)$ is an irreducible unitary representation of $G$. Let $\eta \in \mc H^{-\infty}$ and $\eta \neq 0$. Then
 $$\mc M_{\eta}: (\mc H^*)^{-\infty} \rightarrow \mathbf D^{\prime}(G)$$ is injective.
 \end{thm}
 Proof: We prove this by contradiction. Let $\zeta \in (\mc H^*)^{-\infty}$ such that
 $\zeta \neq 0$ and $\mc M_{\eta}(\zeta)=0$. Then for any $f_1, f_2 \in C_c^{\infty}(G)$, we have
 $$0=\langle \pi(f_1) \eta, \pi^c(f_2) \zeta \rangle.$$
 Let $W_1=\{ \pi(f_1) \eta \mid f_1 \in C_c^{\infty}(G) \}$ and $W_2= \{ \pi^c(f_2) \zeta
 \mid f_2 \in C_c^{\infty}(G) \}$. Then $W_1 \perp W_2$. However $W_1$ and $W_2$ are both
 $G$-invariant subspaces of $\mc H$. Hence either $W_1=\{0\}$ or $W_2=\{ 0 \}$. \\
 \\
 Suppose that $W_1=0$. Then $\pi(f_1) \eta =0$ for all $f_1 \in C_c^{\infty}(G)$. 
 In particular, for any $ v \in (\mc H^*)^{\infty}$, we have
 $$ 0= \langle \pi(f_1) \eta, v \rangle=\int f_1(g) \langle \pi(g) \eta, v \rangle d g=\int f_1(g) \langle \eta, \pi^c(g^{-1}) v \rangle d g $$
 $$=\langle \eta, \int f_1(g) \pi^c(g^{-1}) v d g \rangle.$$
 Notice that the space
 $$ \{ \int f_1(g) \pi^c(g^{-1}) v d g \mid v \in (\mc H^*)^{\infty}, f_1 \in C_c^{\infty}(G) \} $$
 $$= \{ \int f(g) \pi^c(g) v d g \mid v \in \mc (H^*)^{\infty}, f \in C_c^{\infty}(G) \}$$
 is dense in the Frechet space $(\mc H^*)^{\infty}$. So $\eta=0$. This is a contradiction.\\
 \\
 Similarly if $W_2=0$, we have $\zeta=0$. This is also a contradiction. Hence the map 
 $\mc M_{\eta}$ is injective. $\Box$
 \subsection{A counterexample}
 Our generalized matrix coefficients are based on the structure of the Frechet space $\mc H^{\infty}$ and its dual $(\mc H^*)^{-\infty}$. It is possible to choose a subspace  $W \subseteq \mc H^{\infty}$ and its dual $W^*$ to define matrix coefficients in the classical sense. However, it may not be possible to define generalized matrix coefficients for $W^*$. \\
 \\
 Consider $G$ acting on $L^2(G)$ from left. The smooth vectors of $L^2(G)$ are all smooth functions of $G$ whose left Lie derivatives are all in $L^2$. Instead of using the space of smooth vectors, we may use $C_c^{\infty}(G)$  and its dual $\mathbf D^{\prime}(G)$ to define matrix coefficients. In this situation, matrix coefficients are simply certain convolution between $C_c^{\infty}(G)$ and $\mathbf D^{\prime}(G)$. These are well-defined smooth functions. However, we cannot define the generalized matrix coefficients, namely, convolutions between $\mathbf D^{\prime}(G)$ and $\mathbf D^{\prime}(G)$ in the general context, unless $G$ is compact. So the choice of $\mc H^{\infty}$ and $(\mathcal H^*)^{-\infty}$ is essential to guarantee the existence of generalized matrix coefficients. 
 \section{Computation and Application}
 In this last section, we shall address the problem of computing generalized matrix coefficients.
 \subsection{Approximating a distribution by smooth vectors}
 Let $\delta$ be the Dirac delta function.
 Let $\exp: \f g \rightarrow G$ be the exponential map. Let $j(x) \in C^{\infty}_c(\f g)$ with compact support $K$ such that 
 \begin{enumerate}
 \item $\exp$ is a diffeomorphism in a neighborhood of $K$;
 \item $\int j(x) d x =1$.
 \end{enumerate}
  Let $j_n(x)=n^{\dim \f g} j(n x)$. Notice that the support of $j_n(x)$ is $\frac{1}{n} K$ and
  $\int j_n(x) d x=1$. Then we have
  \begin{lem} $j_n(x) \rightarrow \delta_0(x)$ as distributions on $\f g$.
  \end{lem}
  Let $J_n(x)$ be the push forward of $j_n(x)$ by the exponential map, multiplied by the Jacobian. Then we have
  $\int J_n(g) d g=1$ and $J_n(g) \rightarrow \delta_e(g)$ as distributions on $G$.
  \begin{thm} Let $(\pi,\mc H)$ be a unitary representation of $G$. For any $\eta \in \mathcal H^{-\infty}$, $\pi(J_n) \eta \in \mc H^{\infty}$. In addition, $\pi(J_n) \eta \rightarrow \eta$
  in $\mc H^{-\infty}$.
  \end{thm}
  We have obtained an uniform approximation of distributions by smooth vectors. Applying this theorem to the Schr\"odinger representation, we recover the classical result that every tempered distribution on $\mathbb R^n$ can be approximated by functions in the Schwartz space (\cite{fo}). \\
  \\
  Combined with Theorem \ref{5.1}, we have
  \begin{thm} For any $\eta \in \mc H^{-\infty}$ and $\zeta \in (\mc H^*)^{-\infty}$, we have 
  $$\mc M_{\pi(J_n) \eta, \zeta} \rightarrow \mc M_{\eta, \zeta}$$
   in $\mathbf D(G)^{\prime}$. Here $\mc M_{\pi(J_n) \eta, \zeta}$ can be identified with certain smooth functions on $G$.
  \end{thm}
  
\subsection{An example: Fourier series as generalized matrix coefficients}
Let us consider the following classical example, namely Fourier series. We are able to recover some well-known results using generalized matrix coefficients. \\
\\
Let $\mathbb T=\mathbb R/\mathbb Z$. Let $\mathcal H=L^2(\mathbb Z)$. We use $(a)$ to denote functions on $\mathbb Z$. Define a unitary representation of $\mathbb T$ on $L^2(\mathbb Z)$:
$$(\pi(T_t) a)_{n}=(a)_{n} \exp 2 \pi i n t, \qquad (n \in \mathbb Z, (a) \in L^2(\mathbb Z), t \in \mathbb R/\mathbb Z).$$
The Lie algebra $\f t$ can be identified with $\mathbb R X$ with
$$(\pi(X) a)_{n}= 2 \pi i n (a)_n,$$
whenever $\pi(X)(a)$ is in $L^2(\mathbb Z)$. Obviously
$$\mc H^{\infty}=\{ (a) \mid \| (a) \|_k^2=\sum_{n \in \mathbb Z} n^{2k} |(a)_n|^2 < \infty \, \forall \, \, k \geq 0 \}.$$
$\mc H^{\infty}$ is a Frechet space equipped with a set of seminorm $\{\| (a) \|_k \mid k \in \mathbb N \}$.
It is easy to see that 
 $(a) \in \mc H^{\infty}$ if and only if for any $k \in \mathbb N$, there is a $C_k >0$ such that
 $|(a)_n| \leq C_k |n |^{-k}$ for all $n \in \mathbb Z$. In addition, one can use $\{ \sup_n ( |n|^k (a)_n) \mid k \in \mathbb N \}$ as a defining set of seminorms for $\mc H^{\infty}$. \\
 \\
 By Prop \ref{elementary}, $\mc H^{-\infty}= U(\f t) \mathcal H$. Therefore $\mathcal H^{-\infty}$ can be identified with those $(a)$ such that for some $k >0$ and $C_k >0$,
 $$|(a)_n| \leq C_k |n|^k \qquad (\forall \, \, n \in \mathbb Z).$$
 In other words, $(a) \in \mc H^{-\infty}$ if and only if $(a)$ is of at most polynomial growth. \\
 \\
 Let $\mathbf 1=( \ldots, 1,1,1, \ldots)$ be the constant vector. Obviously $\mathbf 1 \in (\mathcal H^*)^{-\infty}$. 
 \begin{thm}
 Given $(a) \in \mc H^{-\infty}$, we have, as distributions on $\mathbb T$,
 $$\mc M_{(a), \mathbf 1}=\sum^* (a)_n \exp 2 \pi i n t.$$
 \end{thm}
 Here $\sum^*{}$ converges in the weak sense, namely, for any $f \in C^{\infty}(\mathbb T)$,
 $$\langle \sum^* (a)_n \exp 2 \pi i n t, f \rangle=\sum \langle (a)_n \exp 2 \pi i n t, f(t) \rangle.$$
 It is not hard to see that $\sum^* (a)_n \exp 2 \pi i n t \in \mathbf D(\mathbb T)^{\prime}$. The proof is straight forward by evaluating both distributions on the test function $f$. 
 \begin{cor}\label{7.1} Let $\{ (b_m) \}$ be a sequence of smooth vectors in $\mathcal H^*$ such that $(b_m) \rightarrow \mathbf 1$ under the weak star topology of $(\mathcal H^*)^{-\infty}$. Then 
 $$\mc M_{(a), (b_m)} \rightarrow \mc M_{(a), \mathbf 1}$$
 in $\mathbf D(\mathbb T)^{\prime}$.
 \end{cor}
 Of course, the convergence here is the weak convergence. We may allow $(b_m)$ to be distributions in $(\mc H^*)^{-\infty}$. We have
 \begin{cor} Let $\{ (b_m) \}$ be a sequence of distributions in $(\mathcal H^*)^{-\infty}$ such that 
 \begin{enumerate}
 \item For each $k$, $(b_m)_k \rightarrow \mathbf (b)_k$, i.e.,  $\{ (b_m) \}$ converges pointwise to a sequence $(b)$;
 \item The sequences $(b_m)$ are uniformly bounded by a sequence $(b_0)$ which is of polynomial growth.
 \end{enumerate}   
 Then 
 $$\mc M_{(a), (b_m)} \rightarrow \mc M_{(a), (b)}$$
 in $\mathbf D(\mathbb T)^{\prime}$.
 \end{cor}
 Proof: (1) and (2) implies that $(b_m) \rightarrow (b)$ in $(\mc H^*)^{\infty}$. Our assertion follows from Theorem \ref{5.1}. $\Box$\\
 \\
In Cor. \ref{7.1}, different choices of $(b_m)$ corresponds to different summation or approximation method for the computation of Fourier series. For example, we may simply choose $(\mathbf 1_m)$ to be the truncation of $\mathbf 1$, namely, $(\mathbf 1_m)_k=1$ if $k \in [-m, m]$ and $(\mathbf 1_m)_k=0$ if $ |k| > m$. Obviously, $(\mathbf 1_m) \rightarrow \mathbf 1$ as distributions. By Cor. \ref{7.1}, we have
\begin{cor} Let $(a) \in \mathcal H^{-\infty}$. Then in $\mathbf D^{\prime}(\mathbb T)$ we have
$$\lim_{m \rightarrow \infty} \sum_{n=-m}^m (a)_n \exp 2 \pi i n t= \sum^* (a)_n \exp 2 \pi i n t.$$
\end{cor}

\end{document}